\documentclass[journal,transmag]{IEEEtran}

\usepackage{amsmath}
\usepackage{graphicx}
\usepackage{amssymb} 
\usepackage{verbatim}
\usepackage{color}
\usepackage{pgfkeys}

\usepackage{amsmath}
\usepackage{ntheorem}
\usepackage{textpos}

\usepackage{amssymb}
\usepackage{bigints}

\usepackage{amsfonts} 
\usepackage{subfigure}

\usepackage{tikz}

\newtheorem{remark}{Remark} 

\newcommand{\naKLPV}{n_{\mathrm{a}_K}}
\newcommand{\nbKLPV}{n_{\mathrm{b}_K}}

\newcommand{\Np}{N_\mathrm{p}}
\newcommand{\Nu}{N_\mathrm{u}}

\newcommand{\aiKLPV}{a^K_i}
\newcommand{\biKLPV}{b^K_i}

\newcommand{\rK}{g}

\newcommand{\matrice}[2]{\left[\hspace*{-.1cm}\begin{array}{#1} #2 \end{array}\hspace*{-.1cm}\right]}

\begin{document}

\title{Direct data-driven control of constrained linear parameter-varying systems:
A hierarchical approach} \vspace{-0cm}


\author{\IEEEauthorblockN{Dario Piga\IEEEauthorrefmark{1},
Simone Formentin\IEEEauthorrefmark{2} and
Alberto Bemporad\IEEEauthorrefmark{1}}
\IEEEauthorblockA{\IEEEauthorrefmark{1}IMT School for Advanced Studies Lucca, Piazza San Francesco 19, 55100 Lucca, Italy}
\IEEEauthorblockA{\IEEEauthorrefmark{2}Politecnico di Milano,  Piazza Leonardo da Vinci 32, 20133 Milano, Italy}    
\thanks{Corresponding author: D. Piga (email: dario.piga@imtlucca.it).}} \vspace{-0cm}

\markboth{IEEE Transactions on Control Systems Technology,~Vol.~??, No.~??}%
{Piga \MakeLowercase{\textit{et al.}}: Direct data-driven control design for LPV systems with constraints}
\IEEEtitleabstractindextext{%
\begin{abstract}
In many nonlinear control problems, the plant can be accurately described by a linear model whose operating point depends on some measurable variables, called scheduling signals. When such a \textit{linear parameter-varying} (LPV) model of the open-loop plant needs to be derived from a set of data, several issues arise in terms of parameterization, estimation, and validation of the model before designing the controller. Moreover, the way modeling errors affect the closed-loop performance is still largely unknown in the LPV context. In this paper, a direct data-driven control method is proposed to design LPV controllers directly from data \emph{without   deriving a model of the plant}. The main idea of the approach is to use a hierarchical control architecture, where the inner controller is designed to match a simple and a-priori specified closed-loop behavior. Then, an outer model predictive controller is synthesized to handle input/output constraints and to enhance the performance of the inner loop. The effectiveness of the approach is illustrated by means of a simulation and an experimental example. Practical implementation issues are also discussed.
\end{abstract}
\begin{IEEEkeywords}
Data-driven control, Linear parameter-varying systems, Constrained control, Model predictive control.  \vspace{-0cm}
\end{IEEEkeywords}}   \vspace{-0cm}

\maketitle

\IEEEdisplaynontitleabstractindextext

\IEEEpeerreviewmaketitle

\section{Introduction}\label{sec:intro}
Linear parameter-varying (LPV) modeling represents an effective tool to describe many nonlinear time-varying systems using linear  \textit{input-output} (IO) maps, wherein changes of an exogenous measurable variable, called scheduling signal, accounts for nonlinear behavior and
time dependency \cite{toth2010modeling}. Using standard robust and gain-scheduling techniques for linear systems, it has been shown that simple and effective controllers can be devised for such complex systems \cite{mohammadpour2012control}.

However, most LPV control design techniques rely on the availability of an accurate physical model of the plant. The effect on the controller of modeling errors between the LPV model and the physical plant is often unpredictable, so that the resulting closed-loop performance might be severely jeopardized. Moreover, even in those applications where gathering data to identify and validate a model of the plant is not costly nor time-consuming, finding a mathematical LPV description of the plant which is good for control design purposes is not an easy task. In fact, when deriving a model of the plant, one always trades off between accuracy and complexity, and, most of the times, is not able to decide a priori how accurate the model should be to achieved a desired closed-loop performance. 


Furthermore, low complexity models of LPV systems are more efficiently derived using input-output model structures \cite{bamieh2002identification,toth2010modeling,piga2015lpv}, which allow to extend Linear Time-Invariant (LTI) prediction-error methods to the  LPV
framework avoiding the curse of dimensionality present in
the identification of state-space LPV models \cite{veve05,felici2007subspace}. On the other hand, most of the control design methods are based on a state-space representation of the system (except some recent works, e.g. \cite{cdc2010,wollnack2013fixed}) and minimal state-space realization of complex IO models is difficult to accomplish \cite{toth2010modeling}.  
 These problems  show that LPV control of nonlinear time-varying models has a great potential, but also suffers from some substantial practical limitations, mostly related to modeling issues.

Recently, a data-driven method has been proposed for directly designing LPV controllers from data, thus avoiding to parametrize, identify and transform an LPV model of the system \cite{formentin2016direct}. This approach sounds appealing and shows many interesting features (e.g., the controller parameters are given by explicit formulas depending on the data points, the mapping with respect to the scheduling signal does not need to be defined a-priori, etc.). However, in some applications this approach cannot be considered as a competitor of other state-of-the-art LPV design techniques, in that signal constraints cannot be taken into account. Furthermore, being a model-reference design method, it requires the desired closed-loop model to be defined, and the choice of an adequate (i.e., practically achievable) reference model without knowing the process dynamics may not be easy. 
These are well-known and open problems in the direct data-driven control literature, both in the LPV and in the LTI framework~\cite{bazanella2011data}.

In this paper, we propose an extension of the data-driven control design method in~\cite{formentin2016direct}. The controller is split into two components, organized in a hierarchical fashion: an inner controller, which accounts for matching a given
simple reference model, and an outer model predictive controller acting as a reference governor~\cite{Bem98}, aiming at enhancing the closed-loop performance and ensure that the constraints are not violated. The main rationale behind this architecture is that the reference model for the inner loop is chosen only to reduce model complexity and uncertainty, but it is decoupled from the desired closed-loop
behavior, which is instead taken care of by the outer part of the controller. Hence, the problem of finding a good reference model becomes less critical than in~\cite{formentin2016direct}
and low-order controller structures can be selected for the identification of the lower-level 
control law from data. Then, the outer model-based controller manipulates the reference signal in such a way that the constraints on input (rate and magnitude) and output are fulfilled
and closed-loop performance increased, without complicating the data-driven design procedure. 
We will show that also the whole control design procedure does not depend on the plant knowledge, according to the direct data-driven philosophy of the method. To the best of the authors' knowledge, this is the first work addressing the problem of handling constraints in direct data-driven control design.  

The effectiveness of the hierarchical control architecture is  illustrated  by means of two examples: ($i$) the simulation case study of~\cite{formentin2016direct}, which best allows us to underline the differences between the proposed method and that of~\cite{formentin2016direct}; ($ii$) an experimental case study concerning the control of an RC circuit with switching load, so as to test the performance of the method when dealing with real-world data. 

The paper is organized as follows. In Section~\ref{sec:prob}, the control problem is formally stated and the additional requirements with respect to~\cite{formentin2016direct} are discussed in detail. The hierarchical architecture of the proposed approach is introduced in Section~\ref{sec:overall_scheme}, while Sections~\ref{sec:inner} and~\ref{sec:outer} discuss the design of the inner and the outer controller, respectively. In the above two sections, methodological details but also practical implementation hints are provided. Finally, the two case studies are illustrated in Section~\ref{sec:case_study}. 

\section{Problem statement}\label{sec:prob}
Let the output signal $y(t) \in \mathbb{R}$, $t\in\mathbb{Z}$, be generated by an unknown \emph{single-input single-output} (SISO) system $\mathcal{G}_{\mathrm{p}}$, driven by the manipulated input
$u(t) \in \mathbb{R}$, a measured exogenous signal $p(t) \in \mathbb{P} \subseteq \mathbb{R}^{n_\mathrm{p}}$, and an unmeasured disturbance $w(t)\in\mathbb{R}^{n_w}$. From now on, we assume $n_\mathrm{p}=1$ to keep the notation simple. 
The system $\mathcal{G}_{\mathrm{p}}$ is assumed to be \emph{bounded-input bounded-output} (BIBO) stable according to the definition in \cite{formentin2016direct}. Assume that
a collection of data $\mathcal{D}_N=\{u(k),y(k),p(k);k \in \mathcal{I}_1^{N}\}$, $\mathcal{I}_1^{N}=\left\{1,\ldots,N\right\}$ generated by the system $\mathcal{G}_{\mathrm{p}}$
is available. 

We aim at synthesizing a controller such that any user-defined (admissible) reference signal can be accurately tracked by the output, without possibly violating the following constraints on inputs and outputs:
\begin{subequations} \label{eq:constraints}
\begin{align}
& u_{\mathrm{min}} \leq u(t) \leq u_{\mathrm{max}}, \ \ \Delta u_{\mathrm{min}} \leq u(t)-u(t-1)  \leq \Delta u_{\mathrm{max}},\\
&  y_{\mathrm{min}} \leq y(t)  \leq y_{\mathrm{max}},\\
&\quad \forall t\in\mathbb{Z},\ t\geq 0.
\end{align}
\end{subequations}
Notice that the (magnitude and rate) constraints on the input are generally imposed by actuator limitations, while the constraints on the output might reflect, for instance, performance specifications or safety conditions. Considering such constraints is therefore of primary importance for many critical engineering applications. 

Rather than attempting at deriving a model of the open-loop plant $\mathcal{G}_{\mathrm{p}}$, we aim at designing a tracking controller directly from the available data set $\mathcal{D}_N$.

\section{A hierarchical approach}\label{sec:overall_scheme}
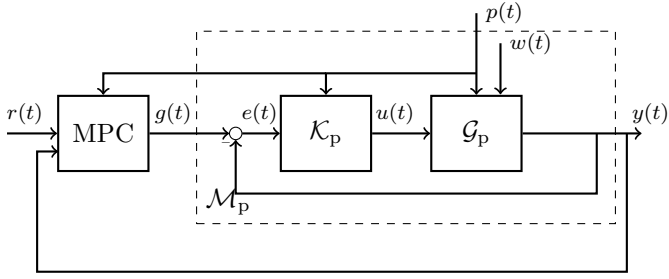
\begin{figure}[!t]
\begin{tikzpicture}[scale=0.8]
\node [draw, rectangle,minimum width=1.2cm,minimum height=1cm,thick] (Gp)  at (4,1.5) {$\mathcal{G}_{\mathrm{p}}$};
\draw[->] [thick] (4,3.5) node[right] {\footnotesize $p(t)$} -- (Gp);
\node [above] (yo)  at (5.2,1.5) {};
\node [draw, rectangle,minimum width=1.2cm,minimum height=1cm, thick] (Kp)  at (1.5,1.5) {$\mathcal{K}_{\mathrm{p}}$};
\draw[<-] [thick] (Kp)--(1.5,2.5)--(4,2.5);
\draw[->] [thick]   (Kp) -- (Gp);
\node [] (plus)  at (6.9,1.5) {};
 \node [draw, circle,scale=0.5] (pluserr)  at (0,1.5) {};
\node [above right] (r)  at (-1.5,1.5) {\footnotesize $g(t)$};
\node [below left,scale=0.5]    at (0,1.5) {$-$};
\draw[->] [thick] (Gp) -- (plus);
\node [right] (ww)  at (4.4,3) {\footnotesize $w(t)$};
\draw[->] [thick] (4.4,3) -- (4.4,2.15);
\node[above] (ynode22) at  (6.9,1.5) {\footnotesize  $y(t)$};
\draw[->] [thick]  (6.5,1.5) -- (6.5,-0.8) -- (-3.3,-0.8) -- (-3.3, 1.2) -- (-2.95,1.2);
 \draw[->] [thick] (pluserr) -- (Kp);
\node [above] (u)  at (2.65,1.5) {\footnotesize $u(t)$};
 \node [above] (e)  at (0.4,1.5) {\footnotesize $e(t)$};
\draw[<-] [thick] (pluserr) -- (0,0.5) -- (6.0,0.5) -- (6.0,1.5);
 \draw[dashed] (-0.65,0) rectangle(6.31,3.2);
\node [above right] (M)  at (-0.65,0) {$\mathcal{M}_{\mathrm{p}}$};
 \node [draw, rectangle,minimum width=1.2cm,minimum height=1cm, thick] (MPC)  at (-2.2,1.5) {$\mathrm{MPC}$};
\draw[<-] [thick] (MPC)--(-2.2,2.5)--(1.5,2.5);
\draw[->] [thick] (MPC) -- (pluserr);
 \node [left] (rout)  at (-3.8,1.5) {};
 \node [above] (rout2)  at (-3.5,1.5) {\footnotesize $  r(t)$};
\draw[->] [thick] (rout)--(MPC);
\end{tikzpicture}  \vspace*{-0.2cm}
\caption{The proposed hierarchical control architecture: the inner controller $\mathcal{K}_{\mathrm{p}}$ provides minimal tracking capabilities for the unconstrained LPV system $\mathcal{G}_{\mathrm{p}}$, whereas the outer MPC controller enhances the performance and guarantees that the constraints are not violated. $\mathcal{K}_{\mathrm{p}}$ is designed from data so that the dynamics of the inner loop is accurately described by $\mathcal{M}_{\mathrm{p}}$.}  \label{fig:hierarchical_outinner}  \vspace*{-0.15cm}
\end{figure}

The proposed control design approach relies on the hierarchical (two degrees of freedom)  architecture illustrated in Fig.~\ref{fig:hierarchical_outinner}, which integrates:
\begin{itemize}
	\item an \emph{inner LPV controller} $\mathcal{K}_{\mathrm{p}}(\theta)$ described by: 
	 \begin{align}\label{eq:K}
 A_K(p,t,q^{-1},\theta)u(t)=B_K(p,t,q^{-1},\theta)(\rK(t)-y(t)),
\end{align}
where
\begin{align} \label{eq:Ak}
 A_{K}(p,t,q^{-1},\theta)= 1+\sum_{i=1}^{\naKLPV}\aiKLPV(p,t,\theta)q^{-i},
\end{align}
\vspace{-3mm}
\begin{align} \label{eq:Bk}
 B_{K}(p,t,q^{-1},\theta)= \sum_{i=0}^{\nbKLPV}\biKLPV(p,t,\theta)q^{-i}.
\end{align}
The dynamical order of the LPV controller $\mathcal{K}_{\mathrm{p}}(\theta)$, defined by the parameters $\naKLPV$ and $\nbKLPV$, is a-priori specified by the user, while  $\aiKLPV(p,t,\theta)$ and $\biKLPV(p,t,\theta)$ are  nonlinear (possibly dynamic) functions of the scheduling variable sequence $p$ and depend on the design parameter vector $\theta$.

The inner controller is designed to achieve a desired LPV (or LTI) closed-loop behavior $\mathcal{M}_{\mathrm{p}}$, a-priori specified by the user and described by the state-space model
\begin{equation} \label{eq:stateMtot}
\begin{array}{rcl}
x_M(t+1) & = & \bar A_M(p,t)x_M(t) + \bar B_M(p,t)\rK(t),  \\
y_{\mathrm{d}}(t) & = & \bar C_M(p,t)x_M(t),
\end{array}
\end{equation}
where $y_{\mathrm{d}}$ denotes the desired closed-loop output for a given reference signal $\rK$.  
The controller parameters $\theta$ achieving the chosen reference model  $\mathcal{M}_{\mathrm{p}}$, as well as the functional dependence on $p$, are estimated directly from the training data set $\mathcal{D}_N$, without first identifying a model for the plant $\mathcal{G}_{\mathrm{p}}$. Such a  data-driven    procedure for LPV control design was originally  introduced in~\cite{formentin2016direct}, and it  will be  reviewed  in Section~\ref{sec:inner}.

\item an \emph{outer LPV model predictive control (MPC) block}, designed based on the desired LPV  closed-loop model $\mathcal{M}_{\mathrm{p}}$. The MPC controller selects, on-line and according to a \emph{receding horizon} strategy, the optimal reference supplied to the inner closed-loop system in order to fulfill the constraints~\eqref{eq:constraints}, thus acting
as a reference governor. Besides constraint fulfillment, the outer MPC allows one
to enhance the performance of the inner closed-loop system modeled 
by $\mathcal{M}_{\mathrm{p}}$. 

\end{itemize}

We will show that the hierarchical structure is an effective choice to solve the direct data-driven \emph{constrained} LPV control problem.
As a matter of fact, on the one hand, the approach in \cite{formentin2016direct} suffers from the drawback that it is difficult to establish whether the selected reference model $\mathcal{M}_{\mathrm{p}}$ is achievable since the plant is unknown. Moreover, there is no way to take the constraints into account. On the other hand, the MPC-based controller alone would need an accurate model of the plant to control, thus a system model should be parameterized, identified and validated. 

By merging the two controllers together in the above hierarchical fashion, one can choose  
a low-demanding (e.g., with slow dynamics and low damping factor) inner closed-loop behavior $\mathcal{M}_{\mathrm{p}}$, which is known to be easily achievable by the inner LPV controller  $\mathcal{K}_{\mathrm{p}}(\theta)$ (for this, only a rough knowledge of the process dynamics is required). The tasks of optimizing the closed-loop performance and fulfilling the input/output constraints are then left to the outer MPC-based controller, which can be designed based on the (known) closed-loop dynamics $\mathcal{M}_{\mathrm{p}}$.

%


\section{Inner controller design}\label{sec:inner}
The main ideas behind the direct data-driven  approach introduced in \cite{formentin2016direct} and employed in this work to design the inner LPV controller $\mathcal{K}_{\mathrm{p}}(\theta)$   are  briefly recalled here for self-consistency of the paper. The design of the outer MPC-based controller is instead discussed in Section~\ref{sec:outer}.

Based on the available training data set $\mathcal{D}_N$, 
the objective is to design a LPV controller $\mathcal{K}_{\mathrm{p}}(\theta)$ achieving a  desired  closed-loop behavior $\mathcal{M}_{\mathrm{p}}$ a-priori specified by the user and described by the state-space equations~\eqref{eq:stateMtot}. 
Unlike~\cite{formentin2016direct}, no specific requirement on the performance of the (inner) closed-loop behaviour $\mathcal{M}_{\mathrm{p}}$ is needed, as the outer MPC will handle the performance requirements. The only assumption that needs to be satisfied by $\mathcal{M}_{\mathrm{p}}$ is that such a behavior is practically achievable. 
Note that this assumption is barely satisfied when a closed-loop $\mathcal{M}_{\mathrm{p}}$ with high performance (e.g., systems exhibiting a high bandwidth and a low overshoot) is chosen. In other words, the chosen LPV controller parametrization might not be flexible enough to directly achieve the desired closed-loop behavior. It is then advisable to impose a low-performance closed-loop behaviour $\mathcal{M}_{\mathrm{p}}$.

\begin{remark}
The above observation can be further clarified by considering a simple LTI example. Consider a model matching problem for a non-minimum phase plant, in which the reference model does not contain the non-minimum phase zeroes of the plant. If the desired bandwidth is high, it is well known that the optimal controller will be likely to destabilize the system in closed-loop \cite{nijmeijer1998approximate}. However, a reference model with a lower bandwidth could still be achieved, as far as the non-minimum phase zeroes are left beyond the desired cut-off frequency.
\end{remark}

In the following, the operator $M(p,t,q^{-1})$ will be used as a shorthand form to indicate the mapping of $g$ to $y_{\mathrm{d}}$ via the reference model $\mathcal{M}_{\mathrm{p}}$. Formally, $M$ is such that $y_{\mathrm{d}}(t)=M(p,t,q^{-1})g(t)$ for all trajectories of $p$ and $g$. Further, we define the left inverse of $M(p,t,q^{-1})$ as the LPV mapping
$M^{\dagger}(p,t,q^{-1})$ that gives $g$ as output when fed by $y_{\mathrm{d}}$, for any trajectory of $p$, i.e., $M^{\dagger}(p,t,q^{-1})M(p,t,q^{-1})=1$.\footnote{For reference maps given in the state-space form~\eqref{eq:stateMtot}, the left inverse $M^{\dagger}(p,t,q^{-1})$ can be computed as indicated in \cite{formentin2016direct}.}

Let $\varepsilon=y_{\mathrm{d}}-y$ be the error between the desired and actual output in response to $g$. According to Fig.~\ref{fig:hierarchical_out2}, we have
\begin{subequations}
\begin{equation}\label{eq:rf}
g(t)\!=\!M^{\dagger}(p,t,q^{-1})y_{\mathrm{d}}(t)\!=\!M^{\dagger}(p,t,q^{-1})(\varepsilon(t)\!+\!y(t)),
\end{equation}
and
\begin{equation}
A_K(p,t,\theta)u(t)= B_K(p,t,\theta)(g(t)-y(t)),
\label{eq:uKp}
\end{equation}
\end{subequations}
$\forall t \in \mathcal{I}_1^N$.
Thus, the controller parameters $\theta$ are computed by minimizing the $2$-norm of the error $\varepsilon$ subject to \eqref{eq:uKp} and \eqref{eq:rf}, \emph{i.e.},
\begin{equation} 
\label{eq:optLPV}
\begin{array}{rl}
\displaystyle \min_{\theta,\varepsilon} &  \sum_{k=1}^N \varepsilon^2(k)\\
\mathrm{s.t.} & A_K(p(k),k,\theta)u(k)\!=\!B_K(p(k),k,\theta)\bigl(M^{\dagger}(p(k),k)\varepsilon(k)\\
&~~+M^{\dagger}(p(k),k)y(k)-y(k)\bigr)
\end{array}
\end{equation}
where $\left\{u(k),y(k),p(k)\right\} \in \mathcal{D}_N$. Notice that problem~\eqref{eq:optLPV} is a purely (non-convex) data-based problem, independent of $\mathcal{G}_{\mathrm{p}}$.
By introducing the residual 
\begin{align}
& \varepsilon_u(\theta,t)=B_K(p(t),t,\theta)M^{\dagger}(p(t),t)\varepsilon(t)=\nonumber\\
&=A_K(p(t),t,\theta)u(t)-B_K(p(t),t,\theta)(M^{\dagger}(p(t),t)y(t)-y(t)),
\label{eq:epsu}
\end{align}
an (approximate) solution of the nonconvex problem~\eqref{eq:optLPV} can be computed, by solving the least-squares problem:
\begin{align} 
\min_{\theta} & \frac{1}{\gamma}\left\|\theta\right\|^2+\frac{1}{N}\sum_{k=1}^N \left|
A_K(p(k),k,\theta)u(k)\right.\nonumber\\
&~~\left.-B_K(p(k),k,\theta)\bigl(M^{\dagger}(p(k),k)y(k)-y(k)\bigr)\right|^2,
\label{eqn:lsestimate}
\end{align}
$\left\{u(k),y(k),p(k)\right\} \in \mathcal{D}_N$, where $\gamma > 0$ is a regularization parameter.  However, since the residuals $\varepsilon_u(\theta,t)$ are not white
 the final estimate of the least-squares problem \eqref{eqn:lsestimate} is not consistent (i.e., the final estimate $\theta$ is not guaranteed to converge to the optimal parameters solving the original problem \eqref{eq:optLPV}) 
and the bias can be not negligible in case of noise $w(t)$ with large variance.
 According to \cite{formentin2016direct}, in order to overcome this problem, the following slight modification of problem \eqref{eqn:lsestimate}, based on instrumental-variables, can be solved instead of  \eqref{eqn:lsestimate}:
\begin{equation} 
\label{eqn:lsestimateIV}
 \min_{\theta,\varepsilon_u} \frac{1}{\gamma}\left\|\theta\right\|^2+\frac{1}{N^2}\left\|\sum_{k=1}^Nz(k)\varepsilon_u(\theta,k)\right\|^2,
\end{equation}
$\left\{u(t),y(t),p(t)\right\} \in \mathcal{D}_N$, where $z(t)$ is the so-called instrument,  chosen by the user so that $z(t)$ is not correlated with the noise $w(t)$. In \cite{formentin2016direct}, it is shown that, in the case $w(t)$ is zero-mean and the output $y(t)$ depends linearly on $w(t)$ (e.g., $w(t)$ is a measurement noise), the final estimate provided by \eqref{eqn:lsestimateIV} converges  to the solution of   problem \eqref{eq:optLPV}.

In the case the controller $p$-dependent coefficient functions $\aiKLPV(p,t,\theta)$ and $\biKLPV(p,t,\theta)$ in \eqref{eq:Ak} and \eqref{eq:Bk} are parametrized as a linear combination of known basis functions of $p$, problems~\eqref{eqn:lsestimate} and
\eqref{eqn:lsestimateIV} are parametric quadratic programming  problems. In the case the dependence of $\aiKLPV(p,t,\theta)$ and $\biKLPV(p,t,\theta)$ on $p$ is not a-priori specified, the dual version of \eqref{eqn:lsestimateIV} can be formulated and the \textit{kernel-based} approaches described in \cite{formentin2016direct} can be used to compute a nonparametric  estimate of the controller coefficients $\aiKLPV(p,t,\theta)$ and $\biKLPV(p,t,\theta)$. 
When Gaussian kernels are used, only the hyper-parameter $\sigma$, representing the width of the  kernels $\displaystyle \kappa(t,j)=e^{\frac{(p(t)-p(j))^2}{\sigma}}$ is specified by the user.

 It is worth mentioning that, it might happen that a part of the control action is a-priori specified, e.g. one may want to include an integrator.  The easiest way to enforce a certain control action $\mathcal{K}^{fixed}_{\mathrm{p}}$ is to process the output data with such a (known) filter, before using the filtered data to identify the remaining part of the controller.
\section{Outer controller design}
\label{sec:outer}
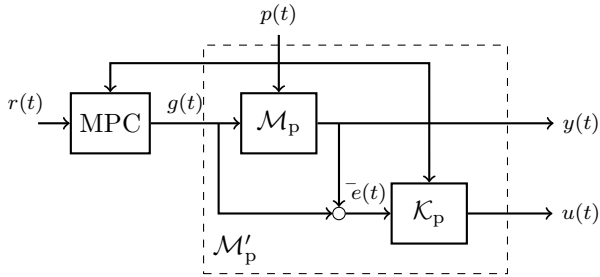
\begin{figure}[!t]
\begin{tikzpicture}[scale=0.8]
\node [draw, rectangle,minimum width=1cm,minimum height=0.8cm,thick] (M)  at (2,1.5) {$\mathcal{M}_{\mathrm{p}}$};
\node [draw, rectangle,minimum width=1cm,minimum height=0.8cm,thick] (MPC)  at (-0.8,1.5) {$\mathrm{MPC}$};
\draw[->][thick] (MPC)--(M);
 \node [left] (rout)  at (-2,1.5) {};
\node [above] (routr)  at (-2.2,1.5) {\footnotesize $r(t)$};
\draw[->][thick] (rout)-- (MPC);
\node [left] (yout)  at (7.5,1.5) {\footnotesize $y(t)$};
\node     at (3.2,0.5) {{\tiny $-$}};
\draw[->] [thick]  (M) -- (yout);
\node [draw, circle,scale=0.5] (minus2)  at (3,0) {};
\node [above]   at (3.5,0) {\footnotesize $e(t)$};
\draw[->] [thick] (3,1.5) -- (minus2);
\node [above] (rMPC)  at (0.45,1.5) {\footnotesize $g(t)$};
\draw[->] [thick]  (1,1.5)-- (1,0) -- (minus2);
\node [draw, rectangle,minimum width=1cm,minimum height=0.8cm,thick] (Kp)  at (4.5,0) {$\mathcal{K}_{\mathrm{p}}$};
\node [left] (uout)  at (7.5,0) {\footnotesize $u(t)$};
\draw[->] [thick]  (minus2)--(Kp);
\draw[->] [thick]  (Kp)--(uout);
\node   (pp) at (2,3.3) {\footnotesize $p(t)$};
\draw [->] [thick] (pp)--(M);
\draw [->] [thick] (2,2.5)--(-0.8,2.5)--(MPC);
\draw [->] [thick] (2,2.5)--(4.5,2.5)--(Kp);
\draw[dashed] (0.75,-1) rectangle(5.8,2.8);
\node[above right] at (0.75,-1) {$\mathcal{M}'_{\mathrm{p}}$};
\end{tikzpicture}  \vspace*{-0.2cm}
\caption{Equivalent single-input two-output LPV model describing the relationship between the MPC output $\rK(t)$ and the plant input and output signals.}   \label{fig:hierarchical_out2}  \vspace*{-0.15cm}
\end{figure}

The outer MPC controller, acting as a reference governor, is designed based on the equivalent single-input two-output model $\mathcal{M}'_{\textrm{p}}$ depicted Fig.~\ref{fig:hierarchical_out2},  where 
 the dynamics of the inner closed-loop system are now described by the (known) model $\mathcal{M}_{\textrm{p}}$. The augmented LPV model $\mathcal{M}'_{\textrm{p}}$ thus describes the relationship between $g(t)$ and $u(t)$, $y(t)$. Within this framework, the role of the inner LPV controller $\mathcal{K}_{\mathrm{p}}(\theta)$ is to transform the behaviour of the unknown plant $\mathcal{G}_{\textrm{p}}$ into that of a known, usually simpler, and  a-priori specified LPV model $\mathcal{M}_{\textrm{p}}$. 

Consider the following, not-necessarily minimal, state-space realization of  $\mathcal{M}'_{\textrm{p}}$
\begin{equation}
\left\{\begin{array}{rcl}
\xi(t+1)&=&A_M(p(t))\xi(t)+B_M(p(t))g(t)\\
\matrice{c}{y(t)\\u(t)}&=&C_M(p(t))\xi(t)+\matrice{c}{0\\D_M(p(t))}g(t),
\end{array}
\right.
\label{eq:SITO}
\end{equation}
where the matrices $A_M(p(t)), B_M(p(t)), C_M(p(t))$ and $D_M(p(t))$ can be easily derived from the  description of the reference model $\mathcal{M}_{\textrm{p}}$ (eq. \eqref{eq:stateMtot}) and the inner controller $\mathcal{K}_{\mathrm{p}}$ (eq. \eqref{eq:K}).
 
Based on the prediction model~\eqref{eq:SITO}, the outer MPC controller is designed both to impose input/output constraints and to possibly improve the tracking quality of the reference signal $r$. As shown in the equivalent scheme of Fig.~\ref{fig:hierarchical_out2}, only the reference model  $\mathcal{M}_{\textrm{p}}$ and the model of the controller  $\mathcal{K}_{\mathrm{p}}(\theta)$ are needed to predict the behaviour of $u(t)$ and $y(t)$. Then, we stress that also in this second step a model of $\mathcal{G}_{\textrm{p}}$ is not required.

\begin{figure*}
\begin{subequations}\label{eq:MPC}
\begin{align}
& \min_{\left\{\rK(t+k|t)\right\}_{k=1}^{\Nu}}Q_y \sum_{k=1}^{\Np} \left(y(t+k|t)-r(t+k)\right)^2+
Q_u \sum_{k=1}^{\Np} \left(u(t+k|t)-u_{\rm ref}(t+k)\right)^2\nonumber\\
&~~~+Q_{\Delta u}\sum_{k=1}^{\Np}\left(\!u(t+k|k)-u(t+k-1|t)\right)^2+ Q_{\rK}\sum_{k=1}^{\Nu}
\left(r(t+k)-\rK(t+k|t)\right)^2+Q_\epsilon \epsilon^2
\end{align}
\begin{align}
\mathrm{s.t. \ } &
	 \xi(t+k+1|t) =  A_M(p(t+k|t))\xi(t+k|t) + B_M(p(t+k|t))\rK(t+k|t), & k=0,\ldots,\Np-1  \\
	& \matrice{c}{y(t+k|t)\\u(t+k|t)} = C_M(p(t+k|t))\xi(t+k|t)+\matrice{c}{0\\D_M(p(t+k|t))}
	\rK(t+k|t),   & k=1,\ldots,\Np  \\
	& -V_y\epsilon+y_{\mathrm{min}}\leq y(t+k|t) \leq  y_{\mathrm{max}}+V_y\epsilon, & k=1,\ldots,\Np \label{eq:MPC:ycons}\\
	& -V_u\epsilon u_{\mathrm{min}}\leq u(t+k|t) \leq  u_{\mathrm{max}}+V_u\epsilon, & k=1,\ldots,\Np  \label{eq:MPC:ucons} \\
	& -V_{\Delta u}\epsilon+\Delta u_{\mathrm{min}}\leq u(t+k|t)-u(t+k-1|t) \leq  \Delta u_{\mathrm{max}}+V_{\Delta u}\epsilon, & k=1,\ldots,\Np \label{eq:MPC:Ducons}\\
	& \rK(t+\Nu+j|t)= \rK(t+\Nu|t), & j=1,\ldots, \Np-\Nu,	\\
	&	\xi(t|t)=\xi(t), \ \ \rK(t)=\rK(t|t). &
\end{align} \vspace*{-5mm}
\end{subequations} \vspace*{-5mm}
\begin{center}
\line(1,0){500}
\vspace{-8mm}
\end{center}
\end{figure*}

The design method is as follows.
By assuming that the state vector $\xi(t)$ of the inner-loop model $\mathcal{M}_{\mathrm{p}}$ 
is fully accessible or, alternatively, estimated from measurements of $u$, $y$ and $p$, for example by means of a linear
time-varying Kalman filter, at each time instant $t$, the reference tracking MPC problem can be formulated, at each time instant $t$, as in~\eqref{eq:MPC}, where $\Np$ and $\Nu$ denote the prediction and control horizon, respectively, $Q_{y}$, $Q_{u}$, $Q_{\Delta u}$, $Q_{\rK}$, $Q_\epsilon$ are nonnegative weights, $u_{\rm ref}$ is a desired input reference (that is typically generated from the output reference $r$ by means of static optimization), $V_y$, $V_u$, $V_{\Delta u}$ are positive vectors that are used to soften the constraints, so that~\eqref{eq:MPC} always admits a solution, that can be computed via 
\textit{Quadratic Programming} (QP).

In the MPC formulation \eqref{eq:MPC},  the following terms are penalized: ($i$) the tracking error between the reference signal $r$ and the  output   $y$; ($ii$)
the tracking error between the input reference signal $u_{\rm ref}$ and the manipulated
variable $u$; ($iii$) the increments of the plant input $u$ (the larger the weight $Q_{\Delta u}$ the less aggressive the control action);  ($iv$) the error between the reference signal $r$ and the MPC output  $g$ and ($v$) the violation of the constraints. From a practical point of view, the goal of the penalty on $g-r$ is to guarantee that the reference signal $g$ of the inner closed-loop system does not differ too much from the reference signal $r$, so as to avoid to excite unmodeled (nonlinear) dynamics.

In case $p(t+k)$ is known at time $t$ for the future $N_p$ steps, we set $p(t+k|t)=p(t+k)$
and call the MPC formulation~\eqref{eq:MPC} \emph{Linear Time-Varying MPC} (LTV-MPC).
In case future values of $p$ are not known, we set $p(t+k|t)\equiv p(t)$ and call the formulation \emph{Linear Parameter-Varying MPC} (LPV-MPC), in which the prediction model is LTI but depends on $p(t)$, and therefore the MPC controller itself is LPV. Alternatively,  the LPV MPC scheme in~\cite{Abbascdc15} can be used to design a robust LPV MPC-based controller. In such an approach, the future values of the scheduling variable are assumed to be uncertain and to vary within a prescribed polytope.

In case both the nominal closed-loop reference model $\mathcal{M}_{\mathrm{p}}$ and the inner controller  $\mathcal{K}_{\mathrm{p}}$  are chosen as LTI models, problem~\eqref{eq:MPC} is a more standard (LTI) MPC problem, that has computational advantages over LTV-MPC and LPV-MPC, in that the QP problem matrices can be precomputed offline, and an explicit MPC approach~\cite{BMDP02a,Bem15} may be viable and reduce the upper control layer to a piecewise affine function.  However, having $\mathcal{M}'_{\mathrm{p}}$ LTI barely happens in practice, in particular when the behaviour of the true plant $\mathcal{G}_\mathrm{p}$ is strongly influenced by the scheduling signal $p$. In this context, even when the selected reference model $\mathcal{M}_{\mathrm{p}}$ is LTI, a parameter-varying controller $\mathcal{K}_\mathrm{p}$ is usually needed to achieve the desired behavior.  

\section{Case studies} \label{sec:case_study}
The effectiveness of the proposed hierarchical control approach 
is shown in this section on two case studies. The first one is the simulation example (concerning the control of a servo positioning system) used in \cite{formentin2016direct} to illustrate the direct data-driven LPV control method. The second case study is an experimental application addressing the control of the output voltage in a RC electric circuit with switching load. These examples show that complex dynamics of  quasi-LPV and switching systems can be dealt with using the approach of the paper. All computations are carried out on  an i7 2.40-GHz Intel core
processor with $4$ GB of RAM running MATLAB R2014b, and the  \textit{Model Predictive Control Toolbox} \cite{BMR15-AdaptiveMPC} is used to design the outer MPC.

\subsection{Simulation case study: the servo positioning system}\label{sec:case_studyDCmotor}

As a first   case study, we consider the  control of a voltage-controlled DC motor with an additional mass mounted on the rotation disc.  In what follows, we show that the hierarchical control structure in Fig.~\ref{fig:hierarchical_outinner} may significantly improve the results  of~\cite{formentin2016direct}, besides allowing us to impose constraints on the input/output signals.

\subsubsection{System description}
The mathematical model of the DC motor, used to simulate the behaviour of the system, is represented by the continuous-time state-space equations
\begin{subequations} \label{eqn:DCmotor}
\begin{align}
\left[\!\begin{array}{c}
  \dot{\theta}(\tau) \\
  \dot{\omega}(\tau) \\
  \dot{I}(\tau)
\end{array}\!\right]\!=\!&
\left[\!
  \begin{array}{ccc}
    0 & 1+\frac{\sin(\theta(\tau))}{\theta(\tau)} & 0 \\
    \frac{mgl}{J}\frac{\sin(\theta(\tau))}{\theta(\tau)} & -\frac{b}{J} & \frac{K}{J} \\
    0 & -\frac{K}{L} & -\frac{R}{L} \\
  \end{array}\!
\right]
\!\!\left[\!\!\begin{array}{c}
  {\theta}(\tau) \\
  {\omega}(\tau) \\
  {I}(\tau)
\end{array}\!\!\right] \nonumber \\
& + \left[
      \begin{array}{ccc}
        0 \ &
        0 \ &
        \frac{1}{L} 
      \end{array}
    \right]^\top V(\tau), \nonumber\\
  y(\tau)= & \left[
            \begin{array}{ccc}
              1 & 0 & 0 \\
            \end{array}
          \right]\left[\begin{array}{c}
  {\theta}(\tau) \\
  {\omega}(\tau) \\
  {I}(\tau)
\end{array}\right], \nonumber
\end{align}
\end{subequations}
where $V(\tau)$ [V] is the control input voltage over the armature, $I(\tau)$ [mA] is the current, $\theta(\tau)$ [rad] is the shaft angle and $\omega(\tau)$ [rad/s] is the angular velocity of the motor. The nomenclature of the parameters characterizing the DC motor is reported in Table \ref{table:DCmotor}, along with their values used to simulate the behaviour
of the motor. The output signal is observed with a sampling time $T_s=10$~ms.

To gather data, the plant  is excited with a discrete-time filtered zero-mean white noise voltage (followed by a zero-order hold block) with Gaussian distribution and standard deviation of $16$ V. The input filter is a first order digital filter with a cutoff frequency of $1.6$ Hz. 
The output measurements are corrupted by an additive white noise $w(\tau)$ with normal distribution and variance such that the Signal-to-Noise Ratio (SNR) is $43$~dB. A second experiment with the same input is also performed to build the instruments $z(k)$ used in \eqref{eqn:lsestimateIV}.


\begin{table}[t!]
\caption{Physical parameters of the DC motor \cite{kulcsar2009lpv}}. \label{table:DCmotor}
\vspace*{-0.6cm}
\begin{center}
\begin{tabular}{l|l|l}
             & Description & Value\\
						 \hline
$R$ &   Motor resistance          &      9.5   $\Omega$         \\
$L$ &   Motor inductance         &      0.84   \!\!$\cdot$\!\! 10$^{-3}$ H         \\
$K$ &   Motor torque constant         &      53.6 \!\!$\cdot$\!\! 10$^{-3}$   Nm/A         \\
$J$ &   Complete disk inertia         &      2.2 \!\!$\cdot$\!\! 10$^{-4}$    Nm$^2$         \\
$b$ &   Friction coefficient         &      6.6 \!\!$\cdot$\!\! 10$^{-5}$  Nms/rad         \\
$M$ &   Additional   mass       &      0.07   kg         \\
$l$ & Mass distance from the  center           &      0.042   m        \\
\hline
\end{tabular}
\end{center} \vspace*{-0.6cm}
\end{table}

\subsubsection{Design of the inner LPV controller $\mathcal{K}_{\mathrm{p}}$} \label{Sec:Ex1inner}
A training data set $\mathcal{D}_N$ with $N=1500$ input/output measurements is used to identify the inner LPV controller $\mathcal{K}_{\mathrm{p}}$  through the procedure discussed in Section~\ref{sec:inner}. The chosen reference model $\mathcal{M}_{\mathrm{p}}$  is described by the  state-space equations:
\begin{equation} \label{eq:exM}
\begin{array}{rcl}
x_M(t+1) & = & 0.99 x_M(t) + 0.01\rK(t) \\
\theta_M(t) & = & x_M(t),
\end{array}
\end{equation}
that is, the desired (inner) closed-loop behaviour $\mathcal{M}_{\mathrm{p}}$ is a simple discrete-time first-order LTI model, with a cutoff frequency of about $6$ Hz.

The  structure for the inner  controller $\mathcal{K}_{\mathrm{p}}$  is given by:
\begin{align*}
u(t)&=\sum_{i=1}^4a^K_i(\Pi(t))u(t-i)+\sum_{j=0}^4b^K_j(\Pi(t))e_{int}(t-j)\\
e_{int}(t)&=e_{int}(t-1)+\left(\rK(t)-y(t)\right),
\end{align*}
where $\displaystyle \Pi(t)=[p(t-1) \ p(t-2) \ p(t-3) \ p(t-4)]^\top$, 
and $p(t)=\theta(t)=y(t)$ (i.e., the output signal measurement  is chosen as scheduling variable). The chosen structure for $\mathcal{K}_{\mathrm{p}}$ is a fourth-order LPV controller   with integral action and  dynamic dependence on the scheduling signal.  
 
 An a-priori parametrization of the coefficient functions $a^K_i(\Pi(t))$ and $b^K_j(\Pi(t))$ is not specified, and Gaussian kernels with width $\sigma=2.4$ are used. The hyper-parameter $\gamma$ in \eqref{eqn:lsestimateIV} is set to $64163$.  The values of   $\gamma$   and the kernel width $\sigma$ are found through cross-validation based on a additional set of $500$ input/output samples. 

\begin{figure}[!t]
\centerline{
\begin{tabular}{c}
\includegraphics[scale=1]{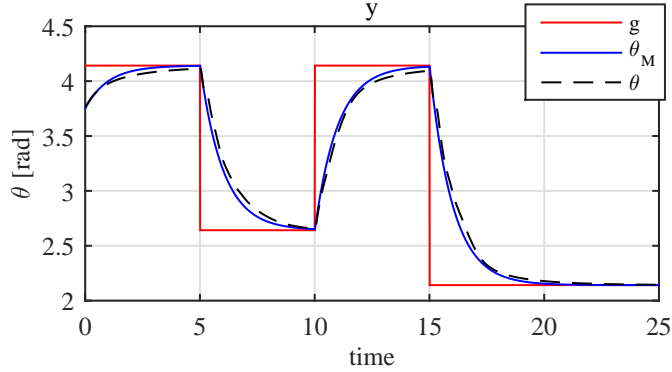}
\end{tabular}}
\vspace*{-0.2cm} \caption{Example 1: inner loop behaviour. Reference signal $\rK(\tau)$ (red), desired step response of the shaft angle $y_d(\tau)$ (solid blue) and actual controlled output $y(\tau)$ (dashed black).} \label{fig:ex1results} \vspace*{-0.2cm} 
\end{figure}

\begin{figure}[!t]
\centerline{
\begin{tabular}{c}
\includegraphics[scale=1]{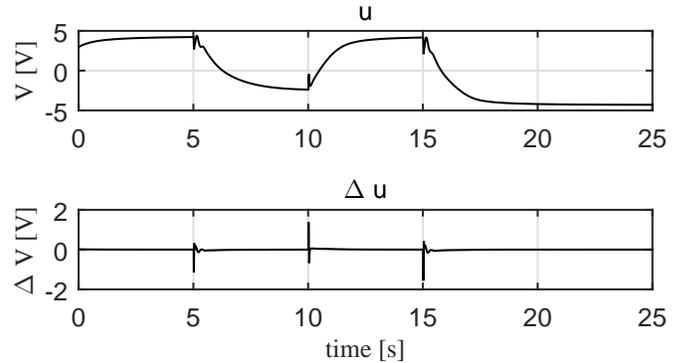}
\end{tabular}} \vspace*{-0.2cm} 
\caption{Example 1: inner loop behaviour. Plant input $V(\tau)$ and input increments $\Delta V(tT_s)=V(tT_s)-V((t-1)T_s)$.} \label{fig:ex1resultsinput} \vspace*{-0.2cm} 
\end{figure}

The performance achieved by the controller $\mathcal{K}_{\mathrm{p}}$ are tested in closed-loop for a piecewise constant reference signal $\rK(t)$. The   response of the (inner) closed-loop system is plotted in Fig.~\ref{fig:ex1results}, and compared  with the output $y_d=\theta_M$ of the desired closed-loop model  $\mathcal{M}_{\mathrm{p}}$  (computed for the same reference excitation). The input voltage $u(t)=V(t)$ provided by the controller $\mathcal{K}_{\mathrm{p}}$ and  applied to the motor  is plotted in Fig.~\ref{fig:ex1resultsinput}.

Results in Fig.~\ref{fig:ex1results} show   a good matching between the actual output $y$ of the closed-loop system and the output $y_d$ of the desired reference model  $\mathcal{M}_{\mathrm{p}}$. However, the closed-loop system exhibits  slow dynamics, with a $10$-$90$\% rise time of about $3.8$~s and a $2$\%-settling time (defined as the time elapsed by the output to enter and remain within a $2$\% error band) of about $4.9$~s.  Unfortunately, due to the limited degrees of the freedom in the controller structure, it has not been possible to achieve  desired reference models  $\mathcal{M}_{\mathrm{p}}$ with   faster  dynamics. A sensitivity analysis with respect to different reference models $\mathcal{M}_{\mathrm{p}}$ is reported in  Table \ref{Tab:ex1MSE}, which  shows the cut-off frequencies of different desired reference models $\mathcal{M}_{\mathrm{p}}$ vs the \textit{mean squares} (MS)
of the differences between the desired closed-loop output $y_d$ and the actual closed-loop output $y$, for the same reference signal plotted in Fig.~\ref{fig:ex1results}. Note that, on the one hand, as the cut-off frequency of the reference model $\mathcal{M}_{\mathrm{p}}$ decreases, the  mismatch between desired and actual closed-loop output decreases, at the price of achieving slower closed-loop dynamics.  On the other hand, for reference models with a cut-off frequency larger than $10$ Hz, the actual closed-loop output $y$ diverges.

\begin{table}
\caption{Cut-off frequency of different reference models $\mathcal{M}_{\mathrm{p}}$ vs mean square (MS) of the difference between desired and actual closed-loop output. The MS is not reported when the achieved closed-loop system is unstable.}
\label{Tab:ex1MSE} \vspace*{-0.2cm}
\begin{tabular}{c|ccccc}
Cut-off & &  &  &\\
frequency [Hz] &  1 & 3 & 6 & 10 & 20   \\ \hline
MS & 0.0001 & 0.0002 & 0.0300& $-$ & $-$  \\
\end{tabular} \vspace*{-0.2cm}
\end{table}

		



\subsubsection{Design of the outer MPC}
Based on the chosen reference model $\mathcal{M}_{\mathrm{p}}$ (which is used to describe the behaviour of the inner closed-loop system) and the designed LPV controller $\mathcal{K}_{\mathrm{p}}$,  an outer MPC is   designed in order to achieve  the following objectives: (i) improve the performance of the inner loop, in terms of rise time and settling time; (ii) enforce the following constraint on the  input voltage rate: $\displaystyle V(tT_s)-V((t-1)T_s) \leq 0.2 V,$ $t=1,2,\ldots$.

The MPC horizons and the weights defining the MPC cost function \eqref{eq:MPC} are tuned through closed-loop simulation, by using  the reference model $\mathcal{M}_{\mathrm{p}}$ to simulate the behaviour of the inner closed-loop system. We stress that this design step is very application-dependent, nevertheless no additional knowledge about the process $\mathcal{G}_{\mathrm{p}}$ is required, being
totally based on the chosen reference closed-loop model $\mathcal{M}_{\mathrm{p}}$. 
The chosen values are equal to $\Np=10$, $\Nu=10$, $Q_y=6.5$, $Q_u=0$, $Q_{\Delta u}=0.1$ and $Q_{\rK}=1$. 

\begin{figure}[!t]
\centerline{
\begin{tabular}{c}
\includegraphics[scale=1]{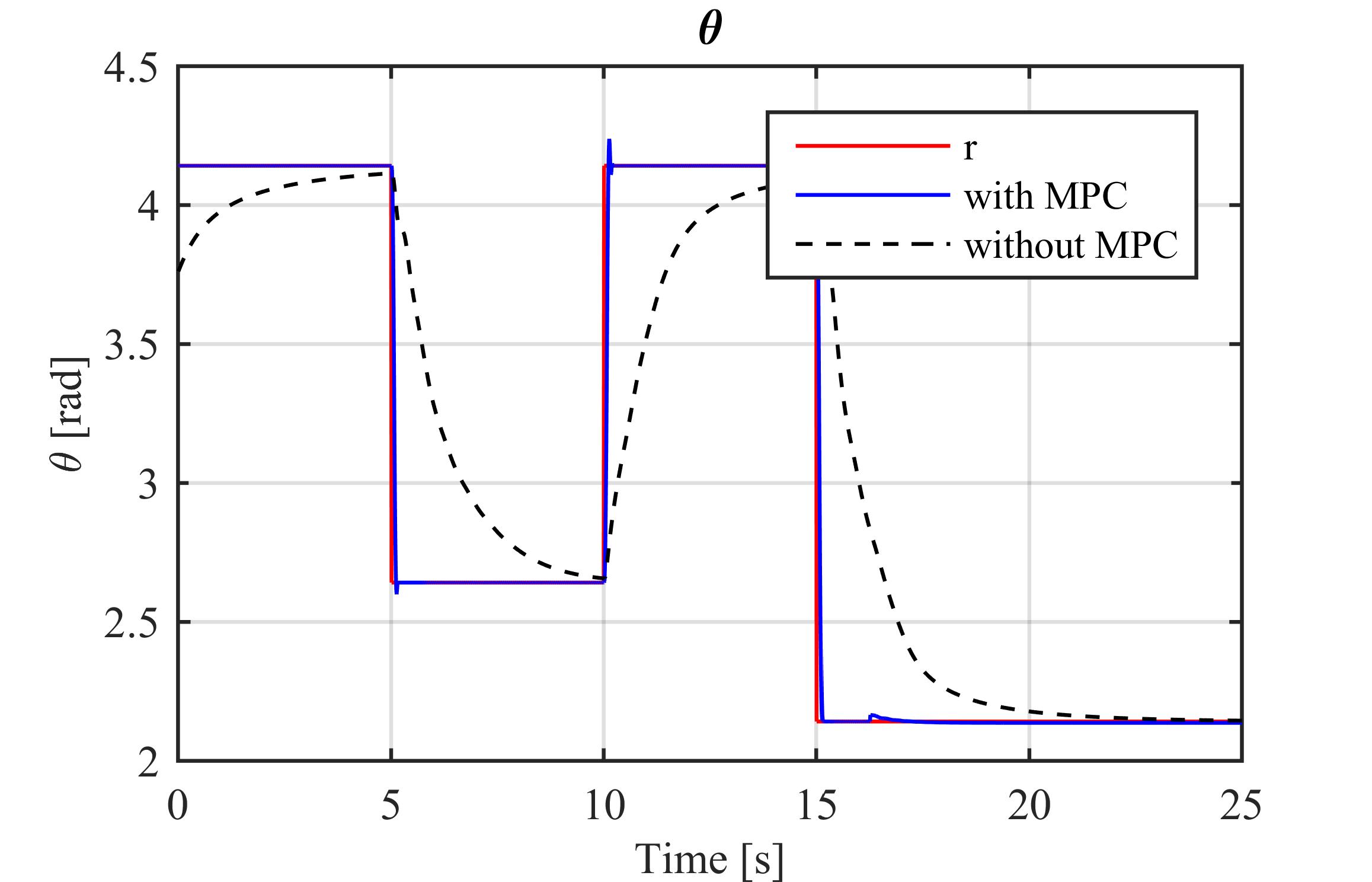}
\end{tabular}}  \vspace*{-0.2cm} \caption{Example 1: closed-loop   behaviour.  Reference signal $r(\tau)$ (red), controlled output $y(\tau)$ (solid blue), and inner-loop output achieved without outer MPC (dashed black).} \label{fig:ex1resultsMPCoutput} 
\end{figure}

\begin{figure}[!t]
\centerline{
\begin{tabular}{c}
\includegraphics[scale=1]{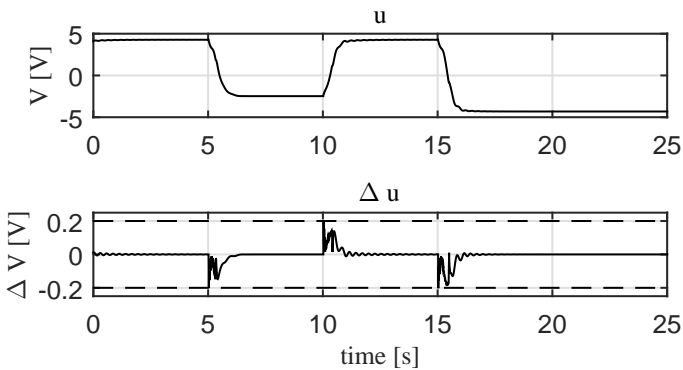}
\end{tabular}}  \vspace*{-0.2cm} \caption{Example 1: closed-loop  behaviour.  Input voltage $V(\tau)$ and input increments $\Delta V(tT_s)=V(tT_s)-V((t-1)T_s)$, constrained between $\pm 0.2$ V (dashed lines).} \label{fig:ex1resultsMPCinputS}
\end{figure}

The response of the closed-loop system  for the same reference signal used in Section \ref{Sec:Ex1inner} is plotted in Fig. \ref{fig:ex1resultsMPCoutput}, while the input voltage applied to the motor is plotted in Fig. \ref{fig:ex1resultsMPCinputS}. For the sake of comparison, the output of the inner  loop achieved without the proposed hierarchical structure is plotted in Fig.  \ref{fig:ex1resultsMPCoutput}. The obtained results show that, although constraints on the variation of the input voltage are enforced,  the hierarchical MPC structure allows us to achieve a faster  reference tracking than the inner-loop system, with a $10$-$90$\% rise time of about $0.7$ s (about $5$ times smaller than the  inner-loop rise time)  and a $2$\%-settling time of about $1.1$ s  (about $4$ times smaller than the inner-loop settling time).
 
The computation time of the MPC layer is 18 ms (including various MATLAB overheads) on the used i7 Intel processor, based on the code generated by the Model Predictive Control Toolbox, which is already in the order of magnitude of the sampling time $T_s=10$~ms. Although computational feasibility is not the main aim of this case study, it is realistic to assume that the controller could be implemented in real-time to control the motor by adopting a fast C implementation of the QP constructor of problem~\eqref{eq:MPC} and QP solver, see the results in~\cite{CBBL15}.

\subsection{Experimental case study: switching RC circuit}\label{sec:case_studyDCmotor}
We address the problem of controlling the output voltage of an RC circuit with switching load.  
 An \textit{Arduino UNO}  board is used for:
(i) measuring the output voltage $V_{\mathrm{out}}$ (namely, the output $y(t)$); (ii) generating the input voltage $V_{\mathrm{in}}$ (namely, the input $u(t)$) applied to the circuit; (iii)  turning on and off the switch (whose driving signal is  the exogenous scheduling signal $p(t)$).

All the computations (including those related to inner and outer control laws) are carried out in MATLAB. The data are transmitted from the Arduino board to MATLAB, and viceversa, via a serial communication at a rate of~9600 baud.


In order to gather the training data set $\mathcal{D}_N$ used to identify the inner control $\mathcal{K}_{\mathrm{p}}$, the following (open-loop) experiment is performed:
\begin{itemize}
\item a  piecewise-constant signal is applied as an input voltage $V_{\mathrm{in}}(t)$ to the electronic circuit;
\item an exogenous piecewise-constant Boolean signal $s(t)$ drives the switch as follows: $s(t)=1$ for Switch  ON, and $s(t)=0$ for Switch  OFF.

\item the  voltage across the capacitor $V_{\mathrm{out}}(t)$ is measured, at a sampling time of $T_s=150$~ms,   with an  \textit{analog-to-digital} (A/D)  converter  available on the Arduino board\footnote{The A/D converters available on the Arduino board used in this experiment have an input rage of $0-5$~V and a resolution of 10~bits.}. A total of~$2000$ samples are acquired, corresponding to a window of $300$~s.  A second measurement of the voltage $V_{\mathrm{out}}(t)$ is taken from another A/D converter to build the instruments. 
\end{itemize}

\begin{figure}[!t]
\centerline{
\begin{tabular}{c}
\includegraphics[scale=1]{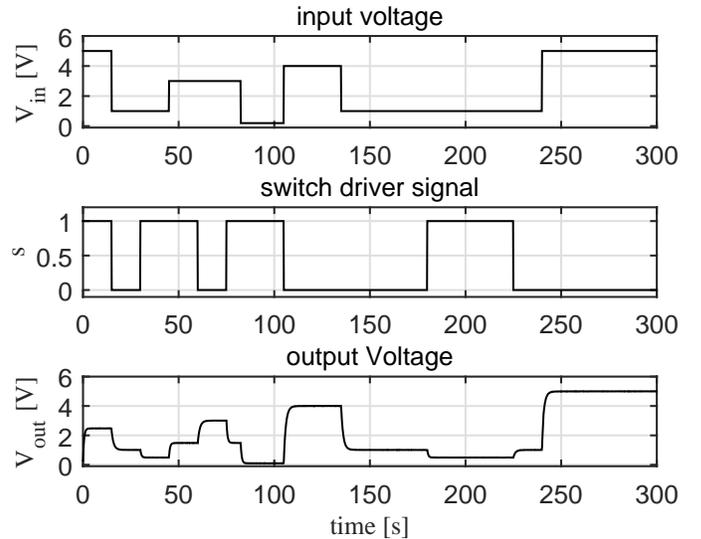}
\end{tabular}}  \vspace*{-0.2cm} \caption{Example 2: open-loop experiment. Input voltage $V_{\mathrm{in}}(\tau)$ (top panel); switch driver signal $s(\tau)$ (middle panel); output voltage  $V_{\mathrm{out}}(\tau)$ (bottom panel).} \label{fig:ex2openloop}
\end{figure}

The signals  $V_{\mathrm{in}}(t)$, $s(t)$ and $V_{\mathrm{out}}(t)$   are plotted in Fig.~\ref{fig:ex2openloop}. 

A new data set with~$500$ samples is also built for tuning the hyper-parameters $\gamma$ and $\sigma$ via cross-validation. 




\subsubsection{Inner LPV controller design}
 The following first-order  LTI model is chosen as a reference model $\mathcal{M}_{\mathrm{p}}$ for the inner loop:
\begin{equation} \label{eq:exM}
\begin{array}{rcl}
x_M(t+1) & = & 0.95 x_M(t) + 0.05\rK(t) \\
\theta_M(t) & = & x_M(t).
\end{array}
\end{equation}
  
A first-order   LPV controller $\mathcal{K}_{\mathrm{p}}$ with an integral action and static dependence on the scheduling variable $p(t)$ is used, i.e., 
\begin{align*}
u(t)&=a^K_1(p(t-1))u(t-1)+\sum_{j=0}^1b^K_j(p(t-1))e_{int}(t-j)\\
e_{int}(t)&=e_{int}(t-1)+\left(\rK(t)-y(t)\right),
\end{align*}
  
The parameters $a^K_1,b^K_1,b^K_2$ defining the LPV controller $\mathcal{K}_{\mathrm{p}}$ are identified through the procedure discussed in Section \ref{sec:inner}. The values of the hyper-parameter  $\gamma$ is $1000$, while kernels  width is $\sigma=1$.



\subsubsection{Design of the outer MPC}
As the Arduino micro-controller can only provide voltage signals within the range $0-5$~V, such a constraint on the signal $u(t)=V_{\mathrm{in}}(t)$ is taken into account while computing the  MPC law for generating $\rK(t)$. Furthermore, the controlled output $y(t)=V_{\mathrm{out}}(t)$ is also constrained to belong to the interval $[0, \ 5]$~V, representing the input range of the A/D converters used in Arduino to measure the voltage  $V_{\mathrm{out}}(t)$.

The following values of the MPC parameters $\Np=3$, $\Nu=3$, $Q_y=0.45$, $Q_u=0$ $Q_{\Delta u}=0$ and $Q_{\rK}=0.1$  are used. These parameters are tuned by means of closed-loop simulations, using  the reference model $\mathcal{M}_{\mathrm{p}}$ as the model of the inner loop. 

\begin{figure}[!t]
\centerline{
\begin{tabular}{c}
\includegraphics[scale=1]{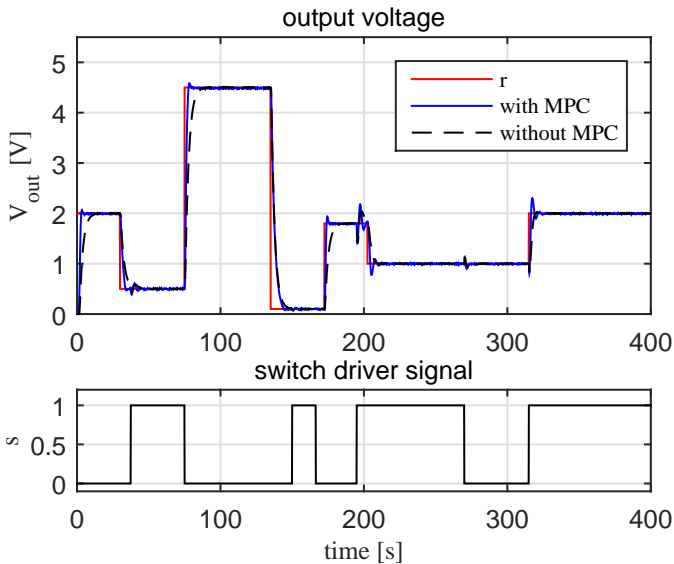}
\end{tabular}}  \vspace*{-0.1cm} \caption{Example 2: closed-loop experiment. Top panel: reference signal (red); controlled output $V_{\mathrm{out}}(\tau)$ (solid blue) and inner-loop output achieved without the outer MPC (dashed black). Bottom panel: switching driver signal $s(\tau)$ during closed-loop experiment.} \label{fig:ex2closedloopVout}
\end{figure}

\begin{figure}[!t]
\centerline{
\begin{tabular}{c}
\includegraphics[scale=1]{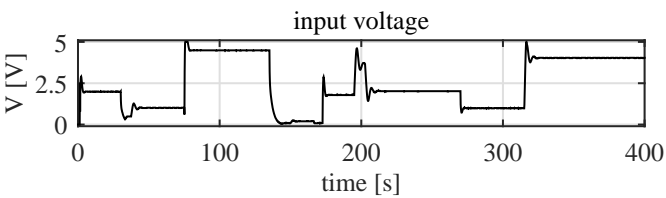}
\end{tabular}}  \vspace*{-0.0cm} \caption{Example 2: closed-loop experiment. Input voltage $V_{\mathrm{in}}(\tau)$.} \label{fig:ex2closedloopu}
\end{figure}

The performance of the designed controllers is then tested by running a closed-loop experiment, with the trajectory of the switching driver signal $s(\tau)$ plotted in Fig.~\ref{fig:ex2closedloopVout} (bottom plot). The obtained controlled output voltage $V_{\mathrm{out}}$ is shown in Fig.~\ref{fig:ex2closedloopVout} (top plot), along with the desired reference signal $r(\tau)$. For the sake of comparison, Fig.~\ref{fig:ex2closedloopVout} also shows the output voltage  $V_{\mathrm{out}}$ achieved by the inner closed-loop system, for the same reference, in the absence of the outer MPC. Notice that such a comparison highlights an evident improvement in terms of raising time for the system with MPC. The trajectories of the input signal $V_{\mathrm{in}}$ is plotted in Fig.~\ref{fig:ex2closedloopu}. The obtained results show that the proposed hierarchical control architecture allows us to efficiently track piecewise constant desired reference voltages in an RC circuit also in the presence of disturbance loads, with faster closed-loop dynamics than the ones achieved by using only the inner LPV controller. Notice that the sudden change of the output load causes only a negligible oscillation on the controlled output voltage   $V_{\mathrm{out}}$  (see Fig.~\ref{fig:ex2closedloopVout} at around $\tau=90$~s and $\tau=320$~s).  
 
The CPU time required to compute the MPC law $\rK(t)$ at each time instant $t$ ranges between $9$~ms and $19$~ms, significantly smaller than the sampling time $T_s=150$~ms. 
 


\section{Conclusions}
In this paper,  a data-driven method to design feedback controllers for LPV systems with constraints is discussed. 
 With respect to the existing works on direct control design available in the literature, constraints on the input and output signals can be accounted for and the choice of the reference model is no longer a critical issue.   To show the effectiveness of the method, we discussed two case studies: the quasi-LPV example in simulation of \cite{formentin2016direct} and an experimental application with a switching RC network. In both the cases, the proposed method shows to be effective and easy to use, and it outperforms the  direct approach of \cite{formentin2016direct}. Future research will deal with: (i) extension of the proposed approach to multivariable systems; (ii) efficient on-line implementation of the outer MPC-based controller; (iii) design of robust controllers to take into account a possible mismatch between the desired and the actual inner closed-loop behaviour. 



\ifCLASSOPTIONcaptionsoff
  \newpage
\fi

\bibliographystyle{IEEEtran}

\bibliography{lpv_AB}

%

\end{document}